\documentclass[12pt,fleqn]{article}

\usepackage{graphicx}
\usepackage{amssymb}

\newtheorem{theorem}{Theorem}[section]
\newtheorem{proposition}[theorem]{Proposition}

\newtheorem{problem}[theorem]{Problem}
\newtheorem{lemma}[theorem]{Lemma}

\newtheorem{corollary}[theorem]{Corollary}
\newtheorem{claim}{Claim}[theorem]
\newtheorem{conjecture}[theorem]{Conjecture}

\begin{document}

\title{Nowhere-zero flows on signed regular graphs}

\vspace{3cm}
      
\author{Michael Schubert\thanks{Fellow of the International Graduate School "Dynamic Intelligent Systems"},  Eckhard Steffen\thanks{
		Paderborn Institute for Advanced Studies in 
		Computer Science and Engineering, 
		Paderborn University,
		Warburger Str. 100,
		33102 Paderborn,
		Germany;			  
		mischub@upb.de; es@upb.de}}      

\date{}

\maketitle

\begin{abstract}
{\small{
We study the flow spectrum ${\cal S}(G)$ and the integer flow spectrum $\overline{{\cal S}}(G)$ of signed $(2t+1)$-regular graphs.
We show that if $r \in  {\cal S}(G)$, then $r = 2+\frac{1}{t}$ or $r \geq 2 + \frac{2}{2t-1}$. 
Furthermore,  
$2 + \frac{1}{t} \in {\cal S}(G)$ if and only if $G$ has a $t$-factor.
If $G$ has a 1-factor, then $3 \in \overline{{\cal S}}(G)$, and for every $t \geq 2$, there is 
a signed $(2t+1)$-regular graph $(H,\sigma)$ with $ 3 \in \overline{{\cal S}}(H)$ and $H$ does not have a 1-factor.

If $G$ $(\not = K_2^3)$ is a cubic graph which has a 1-factor, then $\{3,4\} \subseteq {\cal S}(G) \cap \overline{{\cal S}}(G)$.
Furthermore, the following four statements are equivalent: 
(1) $G$ has a 1-factor. 
(2) $3 \in {\cal S}(G)$.
(3) $3 \in \overline{{\cal S}}(G)$. 
(4) $4 \in \overline{{\cal S}}(G)$.
There are cubic graphs whose integer flow spectrum does not contain 5 or 6, and  
we construct an infinite family of bridgeless cubic graphs with integer flow spectrum $\{3,4,6\}$.

We show that there are signed graphs where the difference between the integer flow number and  the flow number is greater than or
equal to 1, disproving a conjecture of Raspaud and Zhu.

The paper concludes with a proof of Bouchet's 6-flow conjecture for Kotzig-graphs. 
}}
\end{abstract}

\section{Introduction}

Graphs are finite and they may have parallel edges but no loops. The vertex set of $G$ is denoted by $V(G)$, and 
the edge set by $E(G)$. A {\em  signed graph} $(G,\sigma)$ is a graph $G$ and a function $\sigma : E(G) \rightarrow \{ \pm 1 \}$,
which is called a {\em signature} of $G$. The set $N_{\sigma} = \{e : \sigma(e) = -1\}$ is the set of {\em negative edges} of $(G,\sigma)$ and 
$E(G) - N_{\sigma}$ the set of {\em positive edges}.

Let $(G,\sigma)$ be a signed graph.
For $v \in V(G)$, let $E(v)$ be the set of edges which are incident to $v$. 
A {\em switching} at $v$ defines a graph $(G,\sigma')$ with $\sigma'(e) = -\sigma(e)$ for $e \in E(v)$ and $\sigma'(e) = \sigma(e)$ otherwise. We say that
$(G,\sigma)$ and $(G,\sigma^*)$ are {\em equivalent} if they can be obtained from each other by a sequence of switchings. We also say that
$\sigma$ and $\sigma^*$ are equivalent signatures of $G$. 

Every subset $X$ of $E(G)$ defines a signature $\sigma$ of $G$ with $N_{\sigma} = X$. 
If all edges of $G$ are positive, then $N_{\sigma} = \emptyset$ and we say that $\sigma$ is the {\em empty signature}. We 
then write $(G,\emptyset)$ instead of $(G,\sigma)$.

A circuit in $(G,\sigma)$ is {\em balanced}, if it contains an even number of negative edges; otherwise it is {\em unbalanced}. The graph
$(G,\sigma)$ is {\em unbalanced}, if it contains an unbalanced circuit; otherwise $(G,\sigma)$ is {\em balanced}. 
It is well known (see e.g.~\cite{Raspaud_Zhu_2011})
that $(G,\sigma)$ is balanced if and only if it is equivalent to $(G,\emptyset)$. 
Clearly, if  $(G,\sigma)$ and $(G,\sigma')$ are equivalent and $H$ is an eulerian subgraph of $G$, then 
$|N_{\sigma} \cap E(H)|$ and $|N_{\sigma'} \cap E(H)|$ have the same parity. Hence, $H$ is unbalanced in $(G,\sigma)$ if and only if it
is unbalanced in $(G,\sigma')$. 

Let $e \in E(G)$ be an edge, which is incident to the vertices $u$ and $v$. 
We divide $e$ into two half-edges $h_{e}^{u}$ and $h_{e}^{v}$, one incident with $u$ and one incident with $v$. The set of the half-edges of $G$ is denoted $H(G)$. For each half-edge $h \in H(G)$, the corresponding edge in $E(G)$ is denoted by $e_h$. For a vertex $v$, $H(v)$ denotes 
the set of half-edges incident to $v$.
An {\em orientation} of $(G,\sigma)$ is a function $\tau : H(G) \rightarrow \{ \pm 1 \}$ such that $\tau(h_{e}^{u})\tau(h_{e}^{v}) = - \sigma(e)$, for 
each edge $e=uv$.
The function $\tau$ can be interpreted as an assignment of a direction to each edge in the following way. For a  positive edge exactly one half-edge is incoming and the other one is outgoing. For a negative edge either both half-edges are incoming, in which case $e$ is an {\em extroverted} edge, 
or both half-edges are outgoing, in which case $e$ is an {\em introverted} edge. For a vertex $v$ let $H^+(v)$ be the set of outgoing half-edges
and $H^-(v)$ be the set of incoming half-edges which are incident to $v$.

Let $r$, $r'$, $x$ be real numbers with $0\leq r' < r$. We write $ x \equiv r' \pmod r$ if there is 
an integer $t$ such that $x - r' = tr$. 

The {\em boundary} of a function $f : E(G) \rightarrow \mathbb{R}$ is a function $\delta f : V(G) \rightarrow \mathbb{R}$ with 
$\delta f(v) = \sum_{h \in H(v)} \tau(h) f(e_h)$. 
The function $f$ is a ({\em modular}) $r$-{\em flow} on $(G,\sigma)$, if $|f(e)| \in \{x : 1 \leq x \leq r-1\} \cup \{0\}$ for every $e \in E(G)$ and 
	$\delta f(v) = 0$ $(\delta f(v) \equiv 0 \pmod r)$ for every $v \in V(G)$.
The set $\{e : f(e) \not = 0 \}$ is the {\em support} of $f$, and 
$f$ is a {\em nowhere-zero} (modular) $r$-flow on $(G,\sigma)$ if $E(G)$ is the support of $f$.

A signed graph $(G,\sigma)$ is {\em flow-admissible}, if there exists an orientation $\tau$ and a $r \geq 2$ such that  $(G,\sigma)$ has a nowhere-zero 
$r$-flow. 
The {\em circular flow number} of a flow-admissible signed graph $(G,\sigma)$  is 
$$F_c((G,\sigma)) = \inf \{ r : (G, \sigma) \textrm{ admits a nowhere-zero } r \textrm{-flow}  \}.$$
It is known, that $F_c((G,\sigma))$ is always a minimum and that it is a rational number.  
If $(G,\sigma)$ does not admit any nowhere-zero flow, we set $F_c((G,\sigma)) = \infty$. 
If $(G,\sigma)$ and $(G,\sigma')$ are equivalent, then $F_c((G,\sigma)) = F_c((G,\sigma'))$.
Furthermore, $F_c((G,\sigma)) = \infty$
if and only if $(G,\sigma)$ is equivalent to $(G,\sigma')$ with $|N_{\sigma'}| = 1$ or $G$ has a bridge $b$ and
a component of $G-b$ is balanced \cite{Bouchet}. 
If we restrict our studies on integer-valued flows, then $F((G,\sigma))$ denotes the integer flow number of $(G,\sigma)$.
Clearly, $F_c((G,\sigma))  \leq F((G,\sigma))$ and it is known that $F((G,\emptyset)) = \lceil F_c((G,\emptyset)) \rceil$.
Raspaud and Zhu \cite{Raspaud_Zhu_2011} proved that $F((G,\sigma)) \leq 2 \lceil F_c((G,\sigma)) \rceil - 1$, and they
conjectured that $F_c((G,\sigma)) > F((G,\sigma)) - 1$ for every signed graph $(G,\sigma)$. 
In Section \ref{difference} we construct counterexamples to this conjecture. Furthermore, we construct graphs $(G,\sigma)$ with $F_c((G,\sigma)) = 4$ 
and $F((G,\sigma)) = 5$.

Let $G$ be a graph and $X \subseteq E(G)$. Let $\Sigma_X(G)$ be the set of signatures $\sigma$ of $G$, for which $(G,\sigma)$ is flow-admissible
and $N_{\sigma} \subseteq X$. 
We define ${\cal S}_X(G) = \{ r : \textrm{there is a signature } \sigma \in \Sigma_X(G) \textrm{ such that } F_c((G,\sigma)) = r \}$ 
to be the $X$-{\em flow spectrum} of $G$. The $E(G)$-flow spectrum is the {\em flow spectrum} of $G$ and it is denoted by ${\cal S}(G)$.
If we consider integer-valued flows, then $\overline{{\cal S}}_X(G)$ denotes the {\em integer $X$-flow spectrum} of $G$.

The flow number of graphs with empty signature
has been studied intensively in the recent years. One of the most famous conjectures in this context is Tutte's 5-flow conjecture: 

\begin{conjecture} [\cite{Tutte_1954}] \label{5-flow}
Let $G$ be a graph.
If $(G,\emptyset)$ is flow-admissible, then $F((G, \emptyset)) \leq 5$.
\end{conjecture}

Seymour \cite{Seymour_6-flow} proved that $F((G, \emptyset)) \leq 6$ for every bridgeless graph $G$. 
Bouchet formulated the following conjecture for the general case
of signed graphs. 

\begin{conjecture} [\cite{Bouchet}] \label{Bouchet_conj}
Let $(G,\sigma)$ be a signed graph. 
If $(G,\sigma)$ is flow-admissible, then $F((G,\sigma)) \leq 6$.
\end{conjecture}

It is well known that Bouchet's 6-flow and Tutte's 5-flow conjecture are equivalent to their restrictions on cubic graphs. 
The best published approximation to Bouchet's conjecture is proven by  Z\'{y}ka \cite {Zyka_1987}, 
who proved that every flow-admissible graph admits a 
nowhere-zero 30-flow. 

Section \ref{smallest_poss_flow_number} characterizes $(2t+1)$-regular graphs whose flow spectrum contains $2 + \frac{1}{t}$.
Furthermore, if a $(2t+1)$-regular graph has a 1-factor, then its integer flow spectrum contains 3. 
But for every $t \geq 2$, there is 
a signed $(2t+1)$-regular graph $(H,\sigma)$ with integer flow number 3 and $H$ does not have a 1-factor.

One of the earliest results on flows on graphs is Tutte's  characterization of bipartite cubic graphs \cite{Tutte_1949}. His observation that a cubic graph is bipartite if and only if it admits a nowhere-zero 3-flow motivated the following statement. 

\begin{theorem} [\cite{Steffen_2001}] \label{gap_thm} Let $t \geq 1$ be an integer. 
A $(2t+1)$-regular graph $G$ is bipartite if and only if $F_c((G,\emptyset)) = 2 + \frac{1}{t}$. Furthermore,
if $G$ is not bipartite, then $F_c((G,\emptyset)) \geq 2 + \frac{2}{2t-1}$.
\end{theorem}

Section  \ref{gap_signed_graphs} shows
that the situation does not change in the more general case of flow numbers on signed $(2t+1)$-regular graphs. 
We prove that if $r$ is an element of the flow spectrum of a $(2t+1)$-regular graph, 
then $r = 2+\frac{1}{t}$ or $r \geq 2 + \frac{2}{2t-1}$. 

In order to generalize the structural part of Theorem \ref{gap_thm} we will need the following definition:
Let $r \geq 2$ be a real number and $G$ be a graph. A set $X \subseteq E(G)$ is $r$-{\em minimal} if 

1) there is a signature $\sigma$ of $G$ such that $F_c((G,\sigma)) = r$ and $N_{\sigma} = X$, and 

2) $F_c((G,\sigma')) \not = r$ for every signature $\sigma'$ of $G$ with $N_{\sigma'} \subset X$.

In Section \ref{X-spectra} we show that 
a set $X \subseteq E(G)$ is a minimal set such that $G-X$ is bipartite if and only if $X$ is $(2 + \frac{1}{t})$-minimal. 

Since Bouchet's conjecture is equivalent to its restriction on cubic graphs we study 
flows on signed cubic graphs in Section \ref{cubic_graphs}. 
Let $K_2^3$ be the unique cubic graph on two vertices which are connected by three edges. 
We study the relation between 3- and 4-minimal sets and deduce that 
if $G$ has a 1-factor and $G \not = K_2^3$, then $\{3,4\}$ is a subset of its flow spectrum and of its integer flow spectrum.
Furthermore, if $G \not = K_2^3$, then the following four statements are equivalent: 
(1) $G$ has a 1-factor. 
(2) $3 \in {\cal S}(G)$
(3) $3 \in \overline{{\cal S}}(G)$. 
(4) $4 \in \overline{{\cal S}}(G)$.
There are cubic graphs whose integer flow spectrum does not contain 5 or 6, and  
we construct an infinite family of bridgeless cubic graphs with integer flow spectrum $\{3,4,6\}$.

We prove some sharp bounds for the cardinality of
smallest  3-minimal and 4-minimal sets, respectively. If $G$ is not 3-edge-colorable, then these bounds are 
formulated in terms of its resistance and oddness.

 A {\em Kotzig graph} is a cubic graph that has three 
1-factors such that the union of any two of them induces a hamiltonian circuit. 
The paper concludes with a proof of Bouchet's conjecture for Kotzig-graphs.


\section{The difference between the flow number and the integer flow number} \label{difference}

Let $\delta_F = \sup \{F((G,\sigma)) - F_c((G,\sigma)) : (G,\sigma) \textrm{ is flow-admissible}\}$. 
Raspaud and Zhu \cite{Raspaud_Zhu_2011} proved that $F((G,\sigma)) \leq 2 \lceil F_c((G,\sigma)) \rceil - 1$, and
they conjectured that $\delta_F < 1$. We will show that this conjecture is not true.

Let $t \geq 1$ be an integer and $H_t$ be the graph which is obtained from $2t+1$ triangles $T_i$, one vertex $v$ and precisely one vertex of 
each triangle, say $v_i$,  is adjacent to $v$. For $i \in \{1, \dots,2t+1\}$ let $b_i = vv_i$. Clearly, each $b_i$ is a bridge and $H_t$ has no
1-factor. 

\begin{theorem} \label{graph_family}
${\cal S}(H_t) = \{3 + \frac{2}{t}\}$ and $\overline{{\cal S}}(H_t) = \{5\}$, for each integer $t \geq 1$. Furthermore,
$H_t$ has an integer nowhere-zero 5-flow $\phi$ such that $\phi(e) \in \{1, 2, 4\}$,
and a  $(3 + \frac{2}{t})$-flow with $\phi_c$ with $\phi_c(e) \in \{1, 1 + \frac{1}{t}, 2, 2 + \frac{2}{t}\}$ for all $e \in E(H_t)$
\end{theorem}
{\bf Proof.} Let $(H_t,\sigma)$ be flow-admissible. We show that $(H_t,\sigma)$ is equivalent to $(H_t,\sigma^*)$ where precisely the
edge between the two bivalent vertices of each triangle is negative.  
Switch, if necessary, at vertices $v_i$ to obtain an equivalent signature where all bridges are positive. 
Clearly, each triangle is unbalanced. Hence, if three edges of a triangle are negative, then switch at a bivalent vertex such that
precisely one edge of that triangle is negative. Now, if necessary, switch at a bivalent vertex to 
obtain $(H_t,\sigma^*)$. Hence,
$|{\cal S}(H_t)| = |\overline{{\cal S}}(H_t)| = 1$.

\begin{figure}
\centering
	\includegraphics[height=5cm]{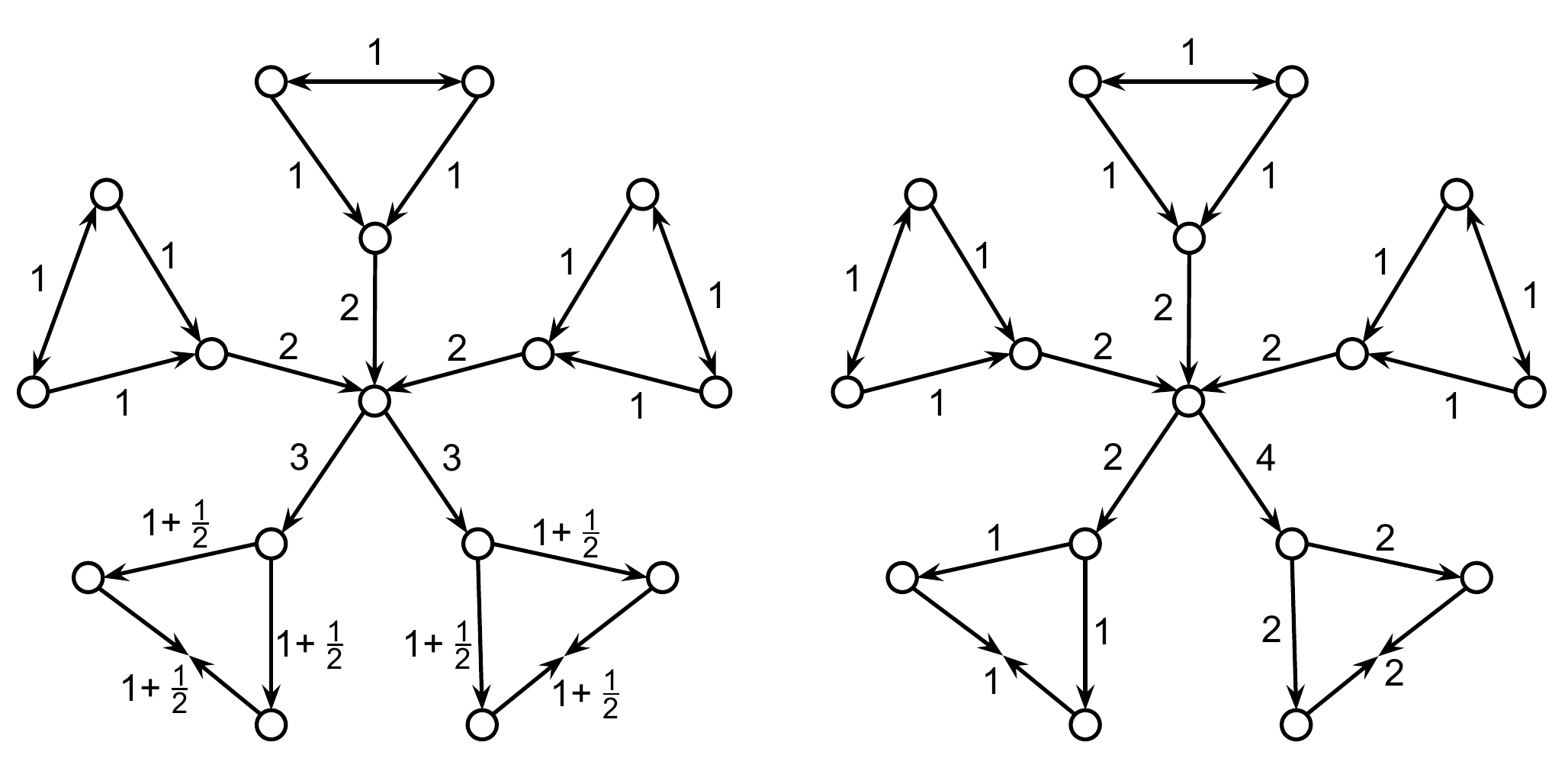}
\caption{$H_2$\label{counterexample}} \end{figure}

We will construct an integer nowhere-zero 5-flow on $(H_t,\sigma^*)$. Define an orientation $\tau$ on $H_t$ as follows.
Let $e_i \in N_{\sigma^*} \cap E(T_i)$ and let
$e_1, \dots, e_{t+1}$ be extroverted and $e_{t+2}, \dots,e_{2t+1}$  be introverted. Orient the positive edges of the triangles
such that the triangles are oriented like a "loop". For $i \in \{1, \dots, t+1\}$, $v$ is the terminal end of  $b_i$, and 
for $j \in \{t+1, \dots, 2t+1\}$, $v$ is the initial end of $b_j$. Let $\phi$ be the integer nowhere-zero 5-flow with $\phi(e)=1$ if 
$e \in \bigcup_{i=1}^{2t}E(T_i)$, $\phi(e) = 2$, if $e \in \{b_1, \dots, b_{2n}\} \cup E(T_{2t+1})$, and $\phi(e) = 4$, if $e = b_{2t+1}$.
Hence, $F((H_t,\sigma^*)) \leq 5$. \\
Let $\psi$ be an integer nowhere-zero flow on $(H_t,\sigma^*)$. Let $E^+(v)$ ($E^-(v)$) be the set
of incoming (outgoing) edges at $v$. Assume that $|E^+(v)| \geq t+1$. Since $\psi$ is an integer flow it follows that $\psi(b_i)$ is even 
for every bridge. Hence, $\sum_{b \in E^+(v)} \psi(b) \geq 2t+2$. Since $|E^-(v)| \leq t$ and $\psi$ is an integer flow it follows that there is a bridge $b$ with 
$\psi(b) \geq 4$. Hence, $F((H_t,\sigma^*)) = 5$, and $\overline{{\cal S}}(H_t) = \{5\}$.

We will construct a nowhere-zero $(3 + \frac{2}{t})$-flow $\phi_c$ on $(H_t,\sigma^*)$. Let $\tau$ be as above and 
$\phi_c(e) = \phi(e)$ if $e \in \bigcup_{i=1}^{t+1}E(T_i)$, $\phi_c(e) = 2$, if $e \in \{b_1, \dots, b_{t+1}\}$,
$\phi_c(e) = 2 + \frac{2}{t}$ if $e \in \{b_{t+2}, \dots,b_{2t+1}\}$, and $\phi_c(e) = 1 + \frac{1}{t}$ if $e \in \bigcup_{i=t+2}^{2t+1}E(T_i)$.
Hence, $F_c((H_t,\sigma^*)) \leq 3 + \frac{2}{t}$. \\
Let $\psi_c$ be a nowhere-zero flow on $(H_t,\sigma^*)$. Assume that $|E^+(v)| \geq t+1$.
Since $\psi_c(b_i) \geq 2$ for every bridge it follows that $\sum_{b \in E^+(v)} \psi_c(b) \geq 2t+2$. 
Hence, there is a bridge $b$ with 
$\psi(b) \geq 2 + \frac{2}{t}$. Therefore, $F((H_t,\sigma^*)) = 3 + \frac{2}{t}$, and ${\cal S}(H_t) = \{3 + \frac{2}{t}\}$.
\hfill $\square$

For $t \geq 2$, the graphs $H_t$ are counterexamples to the conjecture of Raspaud and Zhu. Graph $H_2$ is shown in 
Figure \ref{counterexample}. This graph has also the surprising property that $4 \in {\cal S}(H_2) $, but $ 4 \not \in \overline{{\cal S}}(H_2)$.
As a consequence of Theorem \ref{graph_family} we state:

\begin{corollary} $\delta_F \geq 2$.
\end{corollary}

The statement of Theorem \ref{graph_family} holds also for graphs obtained from $H_t$ by replacing the triangles by (negative) loops. 
Clearly, the argumentation for the lower bounds of the flow numbers in the proof of Theorem \ref{graph_family} works also if replace the 
triangles of $H_t$ by any other unbalanced component and $v$ by a balanced component.

\section{Smallest possible flow numbers of signed $(2t+1)$-regular graphs} \label{smallest_poss_flow_number}

This section characterizes $(2t+1)$-regular graphs whose flow spectrum contains $2 + \frac{1}{t}$.

Let $(G,\sigma)$ be a signed graph that admits a (modular) $r$-flow $\phi$. Let $e$ be an edge of the support of $\phi$. 
If we reverse the orientation of $e$ and replace $\phi(e)$ by $-\phi(e)$ ($r - \phi(e)$) then we obtain another
(modular) $r$-flow $\phi'$ with the same support as $\phi$.
Hence, we can assume that a graph which has a (modular) 
$r$-flow $\phi$ has a (modular) $r$-flow $\phi'$ with the same support and $\phi'(e)  \geq 0$ for all $e \in E(G)$.

\begin{lemma} \label{mod_flow_values_2+1/t}
Let $t \geq 1$ be an integer and $(G,\sigma)$ be a signed $(2t+1)$-regular graph. If $(G,\sigma)$
admits a modular nowhere-zero $(2 + \frac{1}{t})$-flow $\phi$, then 
$|\phi(e)| \in \{1,1 + \frac{1}{t}\}$ for every $e \in E(G)$.
\end{lemma}

{\bf Proof.} Let $\phi$ be a modular nowhere-zero $(2 + \frac{1}{t})$-flow on $(G,\sigma)$.  
Suppose to the contrary that there is an edge $e'$ with $|\phi(e')| \not \in \{1,1 + \frac{1}{t}\}$. 
By the remark above there is    
modular nowhere-zero $(2 + \frac{1}{t})$-flow $\psi$ with 
$\psi(e) > 0$ for every $e \in E(G)$. Since 
$-\phi(e') \not \in \{1,1 + \frac{1}{t}\}$, and $ 2 + \frac{1}{t} - \phi(e') \not \in \{1,1 + \frac{1}{t}\}$
it follows that $\psi(e') \not \in \{1,1 + \frac{1}{t}\}$.

Let $e' \in E(v)$. We can assume that there is an 
orientation of the half-edges of $(G,\sigma)$ such that $|H^+(v)| = 2t+1$ and a modular nowhere-zero $(2+\frac{1}{t})$-flow $\psi'$ such that 
$\psi'(e') \not \in \{1,1 + \frac{1}{t}\}$ and $\psi'(e) > 0$ for every $e \in E(G)$.
It follows that there is a positive integer $k$
such that $\sum_{h \in H^+(v)} \psi'(e_h) = k (2+\frac{1}{t})$. We have  $\sum_{h \in H^+(v)} \psi'(e_h) > 2t+1$ and hence,
$k > t$. On the other side $\sum_{h \in H^+(v)} \psi'(e_h) < 2t+3 + \frac{1}{t}$ and hence, $k < t+1$, a contradiction.
Therefore, $|\phi(e)| \in \{1,1 + \frac{1}{t}\}$ for every $e \in E(G)$.
\hfill $\square$

The next theorem is a generalization of Lemma 3.2 of \cite{XuZhang}.

\begin{theorem} \label{gen_Xu_Zhang}
Let $t \geq 1$ be an integer. A  signed $(2t+1)$-regular graph $(G,\sigma)$ admits a nowhere-zero
$(2+\frac{1}{t})$-flow if and only if $(G,\sigma)$ admits a modular $(2+\frac{1}{t})$-flow and $G$ has a $t$-factor.
\end{theorem}
{\bf Proof.} Let $\phi$ be a nowhere-zero $(2+\frac{1}{t})$-flow on $(G,\sigma)$ with $\phi(e) > 0$ for every $e \in E(G)$. 
Since every $(2 + \frac{1}{t})$-flow is a modular $(2 + \frac{1}{t})$-flow it follows with 
Lemma \ref{mod_flow_values_2+1/t} that $\phi(e) \in \{1,1 + \frac{1}{t}\}$ for every $e \in E(G)$.  
The set of edges with flow value $1 + \frac{1}{t}$ induces a $t$-factor of $G$, and $\phi$ is a modular 
nowhere-zero $(2+\frac{1}{t})$-flow.  

If $(G,\sigma)$ admits a modular nowhere-zero $(2+\frac{1}{t})$-flow, then it follows with Lemma \ref{mod_flow_values_2+1/t}
there is one, say $\phi$, such that $\phi(e) = 1$ for
every $e \in E(G)$. Since $G$ is $(2t+1)$-regular, it follows that every vertex of $G$ is incident to either incoming or outgoing edges, only.
Let $F$ be a $t$-factor of $G$. If we reverse the orientation of the edges of $F$, then the function $\phi'$ with 
$\phi'(e) = 1$ if $e \in E(G) - E(F)$ and $\phi'(e) = 1 + \frac{1}{t}$ if $e \in E(F)$ is a nowhere-zero $(2+\frac{1}{t})$-flow on $(G,\sigma)$.
\hfill $\square$

We will need the following result of Petersen.

\begin{theorem} [\cite{Petersen_1891}] \label{Petersen_2-factor}
Let $k$ be a positive integer and $G$ a $k$-regular graph. If $k$ is even, then $G$ has a 2-factor. 
\end{theorem}

\begin{theorem} \label{always_2+1/t}
Let $t \geq 1$ be an integer. \\
1) A $(2t+1)$-regular graph $G$ has a $t$-factor if and only if $2 + \frac{1}{t} \in {\cal S}(G)$. \\
2) If $G$ is $(2t+1)$-regular and has a 1-factor, then $3 \in \overline{{\cal S}}(G)$. \\
Furthermore, for each $t > 1$ there is a $(2t+1)$-regular
graph $G_t$ which has no 1-factor and  $3 \in \overline{{\cal S}}(G_t)$. 
\end{theorem}
{\bf Proof.} 1) Let $G$ have a $t$-factor. Let $\sigma$ be the signature of $G$ with $N_{\sigma} = E(G)$. 
The function $\phi$ with $\phi(e)=1$ for every edge $e$ 
is a modular nowhere-zero $(2+\frac{1}{t})$-flow on $G$. It follows with Theorem \ref{gen_Xu_Zhang}, that  $2 + \frac{1}{t} \in {\cal S}(G)$.

If  $2 + \frac{1}{t} \in {\cal S}(G)$, then it follows with  Theorem \ref{gen_Xu_Zhang}, that $G$ has a $t$-factor.

2) If $t=1$, then $\phi'$ is an integer 3-flow and it follows from 1) that $3 \in \overline{{\cal S}}(G)$.
Let $t \geq 2$ and $F_1$ be a 1-factor of $G$. By Theorem \ref{Petersen_2-factor}, $G-F_1$
has a 2-factor $F_2$. Hence, $F_1 \cup F_2$ induces a spanning cubic subgraph $H$ of $G$ which has a 1-factor. By 1), 
$H$ has a signature $\sigma$ such that $(H,\sigma)$ has an integer nowhere-zero 3-flow. 
Furthermore, $H'=G - E(H)$ is $(2t - 2)$-regular and hence, $(H',\emptyset)$ has a nowhere-zero 2-flow. 
Thus, $(G,\sigma)$ has an integer nowhere-zero 3-flow. Since a 3-flow is the smallest possible integer flow
on a $(2t+1)$-regular graph it follows that $3 \in \overline{{\cal S}}(G)$.

\begin{figure} \centering
	\includegraphics[height=5cm]{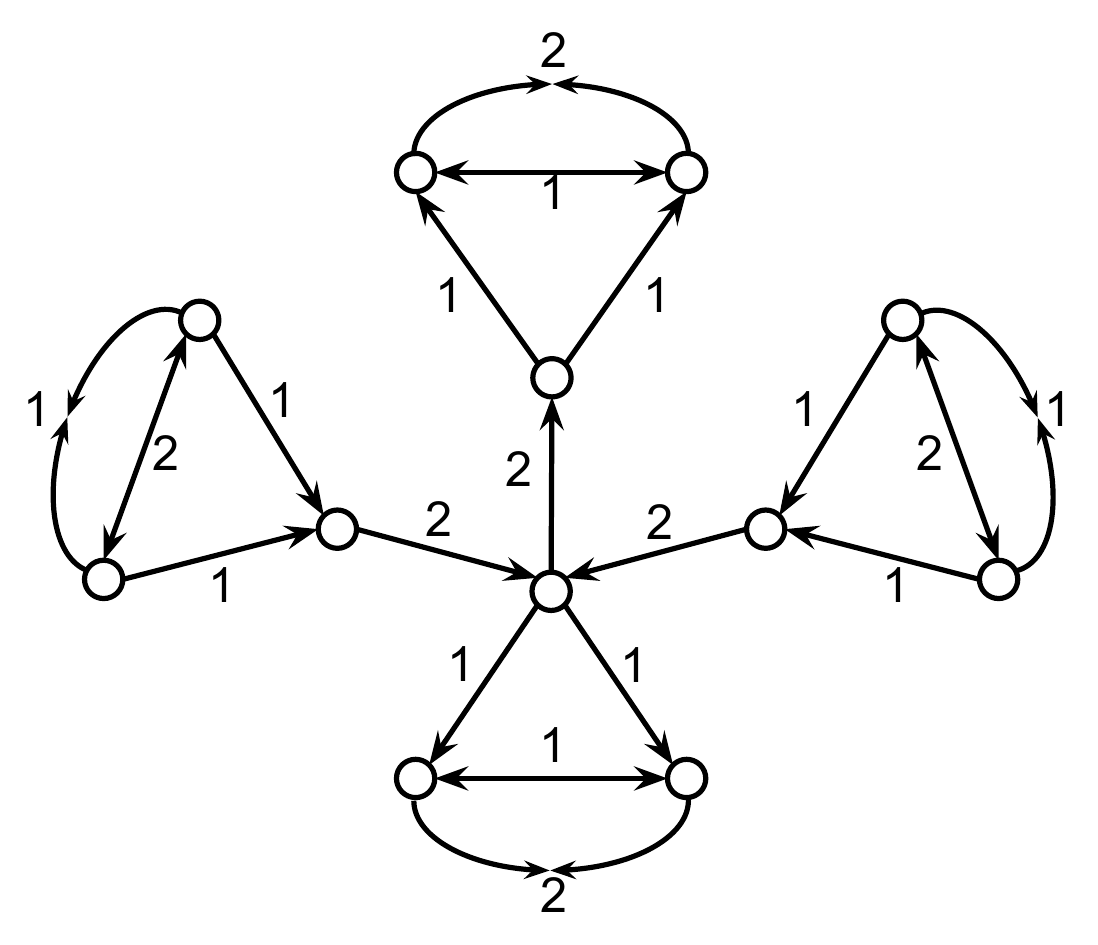}
\caption{A nowhere-zero 3-flow on $(H,\sigma)$ \label{no_1_factor_3_flow}} \end{figure}

It remains to construct the graph $G_t$ for $t>1$.
Let $T$ be a triangle where exactly two vertices are joined by two parallel edges, all other vertices are connected by a simple edge. 
Take four copies $T_1, \dots,T_4$ of $T$ and connect each bivalent vertex of $T_1,T_2,T_3$ with the bivalent vertex $y$ of $T_4$ by an edge. 
Let $H$ be this graph. Let $\sigma$ be the signature of $H$ where $N_\sigma$ is the set of the parallel edges of
$T_1, \dots, T_4$. Graph $(H,\sigma)$ with nowhere-zero 3-flow is shown in 
Figure \ref{no_1_factor_3_flow}.
Since $H$ has a vertex of degree 3, it follows that $F((H,\sigma)) = 3$. 

Let ${K'}_{n,n}$ be the complete bipartite graph on $2n$ vertices where one edge $uv$ is replaced by a path $uxv$. 
For $t > 1$ and for each trivalent vertex $z$ of $H$ take $t-1$ copies of ${K'}_{2t+1,2t+1}$ and identify $z$ and
the bivalent vertices of the copies of ${K'}_{2t+1,2t+1}$. Do the same with $t-2$ copies of  ${K'}_{2t+1,2t+1}$ 
and $y$. The resulting graph $G_t$ is $(2t+1)$-regular. We have $N_\sigma \subseteq E(G_t)$, and  
since 
$F(({K'}_{2t+1,2t+1},\emptyset)) = 3$ it follows that $3 \in \overline{{\cal S}}(G_t)$.

Since $G_t-y$ has more than one odd component it follows that $G_t$ does not have a 1-factor. \hfill $\square$

\section{The set of flow numbers of signed $(2t+1)$-regular graphs} \label{gap_signed_graphs}

This section generalizes the second part of Theorem \ref{gap_thm} to signed graphs. 

\begin{theorem} \label{gap} Let $t \geq 1$ be an integer and $(G,\sigma)$ be a signed $(2t+1)$-regular graph. 
If $F_c((G,\sigma)) = r$, then $r= 2+ \frac{1}{t}$ or $r \geq 2 + \frac{2}{2t-1}$.
\end{theorem}
{\bf Proof.} Let $F_c((G,\sigma)) = r$ and $ r < 2 + \frac{2}{2t-1}$. 
Let $\phi$ be a nowhere-zero $r$-flow on $(G,\sigma)$ and
$\tau$ the corresponding orientation of the half-edges of $(G,\sigma)$ such that $\phi(e) > 0$ for every $e \in E(G)$.  
Let $F_t = \{e : 1 + \frac{1}{2t-1} \leq \phi(e) \leq r-1\}$ and $\overline{F_{t}} = E(G)-F_t$.

Let $v \in V(G)$ and we assume that $|H^+(v)| \geq t+1$
(the statements for the case $|H^+(v)| \leq t$ are proven analogously). 
Let $e \in E(G)$ be an edge which is incident to $v$. 
We will show that $e \in F_t$ if and only if $h_e^v \in H^-(v)$.

 We claim that
$\sum_{h \in H^+(v)} \phi(e_h) < t+1 + \frac{1}{2t-1}$. Suppose to the contrary that 
$\sum_{h \in H^+(v)} \phi(e_h) \geq t+1 + \frac{1}{2t-1}$. Then there is $h' \in H^-(v)$ such that 
$\phi(e_{h'}) \geq \frac{1}{t}(t+1 + \frac{1}{2t-1}) = 1 + \frac{1}{t} + \frac{1}{t(2t-1)} = 1 + \frac{2}{2t-1}> r-1$, a contradiction. 
Thus $|H^+(v)| = |H^-(v)| + 1 = t+1$, and if $e \in F_t$, then $h_e^v \in H^-(v)$. 

For the other direction suppose to the contrary that there is $h' \in H^-(v)$ and $e_{h'} \not \in F_t$. Hence, $\phi(e_{h'}) < 1 + \frac{1}{2t-1}$ and
therefore, $t+1 \leq \sum_{h \in H^+(v)} \phi(e_h) =  \sum_{h \in H^-(v)} \phi(e_h)
	\leq (t-1)(r-1) + \phi(e_{h'})
	< (t-1)(1 + \frac{2}{2t-1}) + 1 + \frac{1}{2t-1}
	= t+1$, a contradiction.
Hence, if $h_e^v \in H^-(v)$, then $e \in F_t$.

Since $|H^+(v)| = |H^-(v)| + 1 = t+1$ it follows that $F_t$ induces a $t$-factor of $G$.
Furthermore, for every vertex $v$ either all half-edges of edges of $F_t$ are incoming and all half-edges of edges of $\overline{F_{t}}$ are outgoing 
or vice versa. Hence, $\phi'$ with $\phi'(e) = 1$ if $e \in \overline{F_{t}}$ and $\phi'(e) = 1 + \frac{1}{t}$ if $e \in F_t$ is a nowhere-zero 
$(2 + \frac{1}{t})$-flow 
on $(G,\sigma)$. Since $2 + \frac{1}{t}$ is the smallest possible flow number for $(2t+1)$-regular graphs it follows that $r =  2 +\frac{1}{t}$.
\hfill $\square$

Theorems \ref{always_2+1/t} and \ref{gap} imply the following corollary.

\begin{corollary}
 Let $t \geq 1$ be an integer and $(G,\sigma)$ be a flow-admissible signed $(2t+1)$-regular graph. 
If $G$ does not have a $t$-factor, then $F_c((G,\sigma)) \geq 2 + \frac{2}{2t-1}$.
\end{corollary}

\section{$r$-minimal sets} \label{X-spectra}

This section studies the structural implications of the existence of a nowhere-zero $(2 +  \frac{1}{t})$-flow on a signed
$(2t+1)$-regular graph. 
Hence, it extends the first part of Theorem \ref{gap_thm} to signed graphs.

\begin{proposition} Let $r \geq 2$ and $G$ be graph. \\
1) The empty set is $r$-minimal if and only if $F_c((G,\emptyset))=r$.\\
2) $F_c((G,\emptyset)) \in {\cal S}_X(G)$ for every $r$-minimal set $X$.
\end{proposition}

Let $G$ be a $(2t+1)$-regular graph and $X \subseteq E(G)$. Let $H$ and $H'$ be two copies of $G-X$. 
For $v \in V(H)$ and $e \in E(H)$, let $v'$ ($e'$) be the corresponding vertex (edge) in $H'$. For each edge $uv \in X$ 
add  edges $uu'$ and $vv'$ (between $H$ and $H'$) to obtain a new $(2t+1)$-regular graph $G_X^2$.
Let $E_X^2$ be the set of the added edges. Let $\sigma$ be a signature on $G$, and $\sigma|_{G-X}$ be the restriction
of $\sigma$ on $G-X$. Let $\sigma_X^2$ be the signature on  $G_X^2$ which  is equal to $\sigma|_{G-X}$
on $H$ and $H'$ and all edges of $E_X^2$ are positive. Note, that 
$|N_{\sigma_X^2}| = 2|N_{\sigma}-X|$. In particular, if $N_{\sigma} \subseteq X$, then $\sigma_X^2$ is the empty
signature.

\begin{lemma} \label{induced_flow}
Let $t \geq 1$ be an integer and $(G,\sigma)$ be a signed $(2t+1)$-regular graph. 
Let $r \geq 2$ and $X \subseteq E(G)$. 
Every nowhere-zero $r$-flow on $(G,\sigma)$ induces a nowhere-zero $r$-flow on $(G_X^2,\sigma_X^2)$.
\end{lemma}
{\bf Proof.} If $\tau$ is an orientation of the half-edges of $(G,\sigma)$, then $\overline{\tau}$ denotes the orientation 
of the half-edges of $(G,\sigma)$ which is obtained from $\tau$ by reversing the orientation of every half-edge. Now,
if $\phi$ is a flow on $(G,\sigma)$ with orientation $\tau$, then $\phi$ is also a flow on $(G,\sigma)$ with orientation $\overline{\tau}$.

If $\phi$ is a nowhere-zero $r$-flow with orientation $\tau$, then define a nowhere-zero $r$-flow on $(G_X^2,\sigma_X^2)$ as follows.
Let $\psi$ be the restriction of $\phi$ on $(H,\sigma)$ with orientation $\tau$ 
and $\psi'$ be the restriction of $\psi$ on $(H',\sigma')$ with orientation $\tau'=\overline{\tau}$.
Extend these orientations to an orientation $\tau_X^2$ on $(G_X^2,\sigma_X^2)$ as follows.
The orientation of the half-edges of $H$ or $H'$ is unchanged, and for an edge $e \in X$ with $e=uv$
orient the edges $uu'$ and $vv'$ of $E_{X}^2$ as follows:
$\tau_X^2(h_{uu'}^{u}) = \tau(h_e^{u})$, 
$\tau_X^2(h_{uu'}^{u'}) = \tau'(h_{e'}^{u'})$, 
$\tau_X^2(h_{vv'}^{v}) = \tau(h_e^{v})$, and
$\tau_X^2(h_{vv'}^{v'}) = \tau'(h_{e'}^{v'})$.
To obtain a nowhere-zero $r$-flow $\phi_X^2$ on $(G_X^2,\sigma_X^2)$ let 
$\psi$ and $\psi'$ be unchanged 
and for $uu', vv' \in E_{X}^2$ which are obtained from
edge $e \in X$ with $e=uv$ let
$\phi_X^2(uu') = \phi_X^2(vv')  = \phi(e)$. \nolinebreak
\hfill $\square$

\begin{theorem} \label{delete_edges_3flow}
Let $t \geq 1$ be an integer and $G$ be a $(2t+1)$-regular graph. 
A set $X \subseteq E(G)$ is $(2 + \frac{1}{t})$-minimal if and only if $G$ has a $t$-factor and $X$ is a minimal set such that $G-X$ is bipartite.
\end{theorem}
{\bf Proof.}  
Let $X$ be $(2 + \frac{1}{t})$-minimal. By definition, there is a signature $\sigma$ of $G$, such that $F((G,\sigma)) = 2 + \frac{1}{t}$, 
$N_{\sigma} = X$ and $F((G,\sigma')) \not = 2 + \frac{1}{t}$ for every signature $\sigma'$ of $G$ 
with $N_{\sigma'} \subset X$. By Theorem \ref{always_2+1/t}, $G$ has a $t$-factor.

Let $\phi$ be a nowhere-zero $(2 + \frac{1}{t})$-flow on $(G,\sigma)$. By Lemma \ref{induced_flow},
$\phi$ induces a nowhere-zero $(2 + \frac{1}{t})$-flow on $(G_X^2,\emptyset)$.
By Theorem \ref{gap_thm}, $G_X^2$ is bipartite and therefore $G-X$ as well.

Now suppose to the contrary that there is a proper subset $X^*$ of $X$ such that $G-X^*$ is bipartite. 

(*) It easily follows that $G_{X^*}^2$ is
bipartite and therefore, $(G_{X^*}^2,\emptyset)$ has a nowhere-zero $(2 + \frac{1}{t})$-flow by Theorem \ref{gap_thm}. 
This $(2 + \frac{1}{t})$-flow can be modified
to a modular nowhere-zero $(2 + \frac{1}{t})$-flow $\phi$ such that $\phi(e)  = 1$ for every edge $e \in E_{X^*}^2$. 
Hence, if we reconstruct $G$ from $G_{X^*}^2$ by keeping the orientation of the half-edges appropriately, we obtain a signature $\sigma^*$
of $G$ with $N_{\sigma^*} \subseteq X^*$ and a modular nowhere-zero $(2 + \frac{1}{t})$-flow $\phi^*$ on $G$. 
By Theorem \ref{gen_Xu_Zhang}, $(G,\sigma^*)$ has a nowhere-zero $(2 + \frac{1}{t})$-flow. 
But $N_{\sigma^*}$ is a proper subset of $X$, contradicting the fact that $X$ is $(2 + \frac{1}{t})$-minimal.   

Let $X$ be a minimal set such that $G-X$ is bipartite. Clearly, if $X \not = \emptyset$, then $|X| \geq 2$. As above (*), we deduce 
that $G$ has a signature $\sigma$ with $N_{\sigma} \subseteq X$
and a nowhere-zero $(2 + \frac{1}{t})$-flow $\phi$ on $(G,\sigma)$. 
Suppose to the contrary that there is an edge $e \in X - N_{\sigma}$. By Lemma \ref{induced_flow}, $\phi$ induces 
a nowhere-zero $(2 + \frac{1}{t})$-flow on $(G_{X-e}^2,\emptyset)$. Hence, $G - (X-e)$ is bipartite, 
contradicting the minimality of $X$. Therefore, $X$ is a $(2 + \frac{1}{t})$-minimal set. \hfill $\square$

An exhaustive survey on sufficient conditions for the existence of factors in regular graphs is given in \cite{Akiyama}.

\section{The flow spectrum of cubic graphs} \label{cubic_graphs}
 
This section studies the flow spectrum and $r$-minimal sets of cubic graphs. 
We will construct some flows on cubic graphs and need the following definition. 
Let $(G,\sigma)$ be a signed graph. For $i \in \{1,2\}$ let $\phi_i$ be a flow on $(G,\sigma)$ with underlying orientation $\tau_i$. 
Note that for each edge $e=uv$ either $\tau_1(h_e^{u}) = \tau_2(h_e^{u})$ and $\tau_1(h_e^{v}) = \tau_2(h_e^{v})$ or
$\tau_1(h_e^{u}) \not = \tau_2(h_e^{u})$ and $\tau_1(h_e^{v}) \not = \tau_2(h_e^{v})$.

The sum $\phi_1+\phi_2$ is the function $\phi$ on $(G,\sigma)$ with orientation $\tau$, where
$\tau = \tau_1|_{\{h| \phi_1(e_h) \geq \phi_2(e_h) \}} \cup \tau_2|_{\{h| \phi_2(e_h) > \phi_1(e_h) \}}$, and
$\phi(e) = \phi_1(e) + \phi_2(e)$ if $e$ has the same direction in $\tau_1$ and $\tau_2$ and 
$\phi(e) = | \phi_1(e) - \phi_2(e) |$, otherwise. Clearly, if $|\phi(e)| \geq 1$ for every edge with $\phi(e) \not = 0$,
then $\phi$ is a flow. 

\begin{theorem} \label{facts_4_minimal_bipartite}
Let $G$ be a cubic graph. 
If $G$ is bipartite, then it has signature $\sigma$ with $|N_\sigma| = 2$ and $F_c((G,\sigma)) = F((G,\sigma)) = 4$. 
\end{theorem}
{\bf Proof.}  By Theorem \ref{gap_thm}, $(G,\emptyset)$ has a nowhere-zero 3-flow $\phi$. 
We can assume that $\phi(e) \geq 1$ for
each $e \in E(G)$. The edges with flow value 1 induce a 2-factor $F$ of $G$. Let $C$ be a circuit of $F$. Any two adjacent edges of $C$ are 
both oriented towards the vertex they share or both away from it. Hence, $\phi$ can be modified to a nowhere-zero 3-flow $\phi'$,
where the circuits of $F$ are directed circuits and the flow values on the edges are 1 and -1, alternately.

\begin{claim} \label{Claim 1} 
There is a circuit $C$ of $F$ and a path $P$ of three consecutive edges $e_1,e_2,e_3$ of $C$ such that
$G-\{e_1, e_3\}$ is connected. Furthermore, $G$ contains a circuit $D$ such that $e_1 \in E(D)$ and $e_3 \not \in E(D)$.
\end{claim}
{\bf Proof.} Every circuit of $F$ has length at least 4. 
Let $C$ be a circuit of $F$, and $e_1= v_1v_2$, $e_2 = v_2v_3$, $e_3 = v_3v_4$ and $e_4 = v_4v_5$ ($v_1=v_5$ is not excluded) 
four consecutive edges
in $C$. Furthermore, $v_2$ has a neighbor $x \not \in \{v_1,v_3,v_4\}$. If $\{e_1,e_3\}$ is a 2-edge-cut in $G$, then choose $P'$ with edges $e_2,e_3,e_4$. Suppose to the contrary that
$\{e_2,e_4\}$ is an edge-cut. It follows that $e_2$ is simple and $v_2x$ is a bridge, contradicting the fact that $G$ is bridgeless. 
Hence, there is a path as claimed.
Furthermore, $G - e_3$ is bridgeless and hence, $G - e_3$ contains a circuit $D$ with $e_1 \in E(D)$. 
$\square$ 

Let $P$ be the path in $C$ with three consecutive edges $e_1,e_2,e_3$ ($e_i=v_iv_{i+1}$) such that $G-\{e_1,e_3\}$ is connected.  
We can assume that $e_2$ is directed from $v_2$ to $v_3$ and that $\phi'(e_2) = 1$. Let $\tau$
be the orientation of $H(G)$ which is obtained from the underlying orientation for $\phi'$ by reversing the 
orientation of $h_{e_1}^{v_1}$ and $h_{e_3}^{v_4}$.
Hence, we obtain a signature $\sigma$ of $G$ with $N_\sigma = \{e_1,e_3\}$, where $e_1$ is extroverted and $e_3$ introverted.
Consider $\phi'$  on $(G,\sigma)$, then
$\delta\phi'(v_1) = -2$, $\delta\phi'(v_4) = 2$, and $\delta \phi'(v) = 0$ if $v \in V(G) - \{v_1,v_4\}$.
Let $E(\overline{P}) = E(C) - E(P)$.  
The
function $\psi :  E(G) \rightarrow \{1,2,3\}$ with $\psi(e) = \phi'(e)$ if $e \in E(G) - E(\overline{P})$ and $\psi(e) = \phi'(e) + 2$ if $e \in E(\overline{P})$
is a nowhere-zero 4-flow on $(G,\sigma)$. Since $\psi$ is an integer flow it follows that $F((G,\sigma)) \leq 4$.

By Claim \ref{Claim 1} and Theorem \ref{gap} it follows that $F_c((G,\sigma)) \geq 4$ and hence, 
$F_c((G,\sigma)) = F((G,\sigma))=4$.\hfill $\square$

We will use the strict form of Petersen's Theorem on 1-factors in cubic graphs. 

\begin{theorem} [\cite{Petersen_1891}] \label{Petersen}
Let $G$ be a bridgeless cubic graph. For every $e \in E(G)$ there is a 1-factor of $G$ that contains $e$. 
\end{theorem}

The minimum number of odd circuits of a 2-factor of a cubic graph $G$ is the {\em oddness} of $G$ 
and it is denoted by $\omega(G)$. A 2-factor that has precisely $\omega(G)$ odd circuits is a minimum 
2-factor of $G$.

\begin{theorem} \label{facts_4_minimal_class2}
Every non-bipartite cubic graph $G$ with 1-factor has a signature $\sigma$ such that $|N_\sigma| = \omega(G)$
and $F_c((G,\sigma)) = F((G,\sigma)) = 4$.
Furthermore, if $G$ is bridgeless, then for every 3-minimal set $X_3$ there is a 4-minimal set $X_4$ with
$X_4 \subset X_3$ and $|X_4| \leq \omega(G) < |X_3|$.
\end{theorem}
{\bf Proof.}   If $G$ is 3-edge-colorable, then the empty set  is 4-minimal and the statements follow. 

Thus, we assume that $G$ is not 3-edge-colorable in the following.
Let $\omega(G) = 2n$, $F_2$ be a 2-factor with odd circuits $C_1, \dots, C_{2n}$, and $F_1$ be the complementary 1-factor. 
If $G$ has a bridge, then there is one, say $b$, such that one component of $G-b$ is bridgeless. Such components will be called
end-components. 

We first show that there is a signature $\sigma$ with $|N_{\sigma}| = \omega(G)$ and $F_c((G,\sigma)) = F((G,\sigma)) = 4$. 
For $i \in \{1, \dots, 2n\}$ choose $f_i \in E(C_i)$ with the following restrictions if
$G$ has bridges or if an odd circuit of $F_2$ contains a multi-edge. 

(1) One of the odd circuits of $F_2$, say $C_k$, has two vertices which are connected by two edges in $G$.
Choose $f_k$ to be one of these two edges, and for $i \not = k$ choose $f_i \in E(C_i)$ arbitrarily. 

(2) All edges of the odd circuits of $F_2$ are simple in $G$ and $G$ has an end-component $K$ such that 
the bivalent vertex $x$ is contained in an odd circuit $C_k$  of $F_2$. Let $x_1$ and $x_2$ be the 
two neighbors of $x$ in $K$. Then choose $f_k$ to be an edge of $C_k$ which is incident to $x_1$ and different from $xx_1$.
For $i \not = k$ choose $f_i \in E(C_i)$ arbitrarily. 

In all other cases choose $f_i \in E(C_i)$ arbitrarily.

Subdivide $f_i$ by a vertex $u_i$ and add edges $e_k = u_{2k-1}u_{2k}$, for $k \in \{ 1,\dots,n\}$. The resulting graph $G'$ is cubic. 
The set $F'_1 = F_1 \cup \{e_k \textrm{ : } k = 1, \dots,n  \}$  is a 1-factor and $F'_2 = E(G') \setminus F'_1$ is an even 2-factor of $G'$.
The odd circuits $C_i$ of $F_2$ are transformed into even circuits $C_i'$ of $F_2'$. Let $f_i'$ and $f_i''$ be the two edges of $C_i'$
which are incident to $u_i$.  
Let $c$ be a proper 3-edge-coloring of $G'$ that colors the edges of $F'_1$ with color $1$, and the edges of the circuits of $F'_2$ with
colors $2$ and $3$. Assume that $f_i' \in c^{-1}(2)$. 
Let $\phi_1$ be a nowhere-zero 2-flow on $G'[c^{-1}(1) \cup c^{-1}(2)]$. Let $\phi_2$ be a nowhere-zero 2-flow on 
$G[c^{-1}(2) \cup c^{-1}(3)]$,
with the additional property that the orientation of $f_i'$ is different in $\phi_1$ and $\phi_2$. Flow $\phi_2$ exists, since 
for every $i \in \{1, \dots,2n\}$ there is precisely one edge $f_i'$ in $C_i'$. 

Now, $2\phi_1 + \phi_2$ is a nowhere-zero 4-flow $\psi$ on $G'$ with the additional property that for each $i \in \{1, \dots,2n\}$ the vertex $u_i$
is either the terminal vertex of $f_i'$ and of $f_i''$ or it is the initial vertex of both of these edges. Furthermore, $\psi(f_i') = \psi(f_i'') = 1$.

For each $i \in \{1, \dots,n\}$, remove edge $e_i$, suppress vertices $u_{2i-1},u_{2i}$ and consider $f_i'$ and $f_i''$ as two half-edges 
of $f_i$ to construct a signature $\sigma$ of $G$ with  $N_{\sigma} = \{f_1, \dots, f_{2n}\}$. Furthermore, $\psi$ induces an integer
nowhere-zero 4-flow on $(G,\sigma)$. Hence, $F((G,\sigma)) \leq 4$.

It remains to show that $F_c((G,\sigma)) \geq 4$. Suppose to the contrary that $F_c((G,\sigma)) < 4$. Then,
$F_c((G^2_{N_{\sigma}})) < 4$ and hence, $G - N_{\sigma}$ is bipartite by Theorem \ref{gap_thm}. 

(1') If there is $k \in \{1, \dots ,2n\}$ such that $C_k$ has two vertices which are connected by two edges in $G$,
then it follows from the construction of $\psi$, that $G - N_{\sigma}$
contains an odd circuit, a contradiction. Hence,  $F_c((G,\sigma)) = F((G,\sigma)) = 4$.

We may assume that all edges of the odd circuits of $F_2$ are simple in $G$. 

(2') If $G$ has an end-component, then there is one, say $K$, which is respected in the construction of $\psi$. 
Let $F_2[E(K)]$ be the 2-factor of $K$ which is a subgraph of $F_2$. 
It follows with Theorem \ref{Petersen} (suppress the bivalent vertex), that there is a 2-factor $F'$ of $K$ that
does not contain $f_k$. Since $K$ has odd order and $F$ contains at least as many odd circuits as $F_2[E(K)]$ 
it follows that $F'$ contains an odd circuit that does not contain any edge of $N_{\sigma}$.
Thus, $G-N_{\sigma}$ contains an odd circuit, a contradiction. Hence,  $F_c((G,\sigma)) = F((G,\sigma)) = 4$.

It remains to consider the case when $G$ is bridgeless. Let $X_3$ be a 3-minimal set. 
Since $X_3$ contains an edge of every odd circuit of $F_2$ it follows that $|X_3| \geq \omega(G)$. 
For $i \in \{1, \dots,2n\}$ let $f_i \in X_3 \cap E(C_i)$. Let $\sigma'$ be the signature
on $G$ with $N_{\sigma'} = \{f_1, \dots,f_{2n}\}$ and 
construct an integer nowhere-zero 4-flow on $(G,\sigma')$ as above. Hence, 
$X_3$ contains a 4-minimal set $X_4$ with $|X_4| \leq \omega(G)$.

We will show that $|X_3| > \omega(G)$.
Suppose to the contrary that $|X_3| = \omega(G)$. Then $\{f_i\} = E(C_i) \cap X_3$ for each $i \in \{1, \dots,2n\}$. 
By Theorem \ref{Petersen} there is a 1-factor of $G$ that contains $f_1$. The complementary 2-factor $F'$ has at least
$\omega(G)$ odd circuits. Thus, there is an odd circuit of $F'$ that does not contain an edge of $X_3$ which implies that 
$G-X_3$ is not bipartite, a contradiction. Thus, $|X_3| > \omega(G)$, $X_4 \subset X_3$, and $F_c((G,\sigma')) = F((G,\sigma')) = 4$. 
\hfill $\square$

Every signature of $K_2^3$ is either equivalent to a signature with no negative edges or to a signature with precisely one negative edge. 
Hence, ${\cal S}(K_2^3) = \overline{{\cal S}}(K_2^3) = \{3\}$. 
The following statement follows with Theorems \ref{always_2+1/t},
\ref{facts_4_minimal_bipartite}, and \ref{facts_4_minimal_class2}. 

\begin{theorem}\label{main_theorem}
Let $G$ be a cubic graph which has a 1-factor. If $G \not = K_2^3$, then $\{3,4\} \subseteq {\cal S}(G) \cap \overline{{\cal S}}(G)$.
\end{theorem}

\begin{theorem}\label{main_theorem_2}
Let $G$ be a cubic graph. If $G \not = K_2^3$, then the following statements are equivalent. \\
1) $G$ has a 1-factor.\\
2) $3 \in {\cal S}(G)$. \\
3) $3 \in \overline{{\cal S}}(G)$.\\
4) $4 \in \overline{{\cal S}}(G)$.
\end{theorem}
{\bf Proof.} 
Statement 1) implies 2) by Theorem \ref{main_theorem}. By Lemma \ref{mod_flow_values_2+1/t} it follows that every circular nowhere-zero 3-flow on 
$G$ is an integer nowhere-zero 3-flow. Hence, statement 2) implies statement 3), which implies statement 1 by Theorem \ref{always_2+1/t}. 

We show the equivalence of 
statements 1) and 4). If $G$ has a 1-factor, then $4 \in \overline{{\cal S}}(G)$ by Theorems \ref{facts_4_minimal_bipartite}
and \ref{facts_4_minimal_class2}. If $4 \in \overline{{\cal S}}(G)$, then $G$ has an integer nowhere-zero 4-flow $\phi$ with
$\phi(e) > 0$ for each edge $e$. It is easy to see that $F = \{e : \phi(e) = 2\}$ is a 1-factor of $G$.\hfill $\square$ 

Theorem \ref{main_theorem} says that if $G$ has a 1-factor, then $4 \in {\cal S}(G)$. However, the other direction is not true.

\begin{figure}
\centering
	\includegraphics[height=5cm]{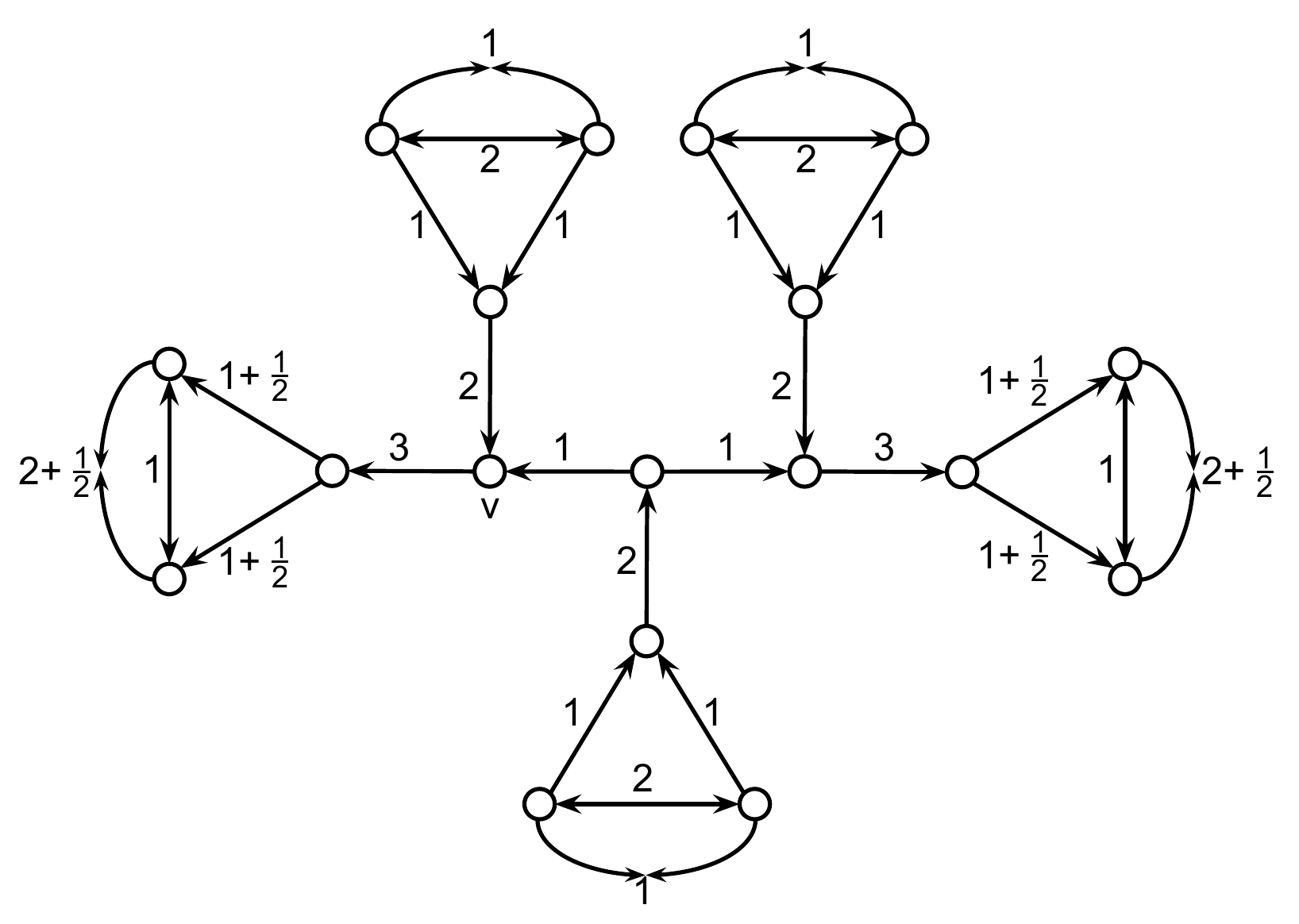}
\caption{A graph $H$ which has no 1-factor and $4 \in {\cal S}(H)$ \label{best}} \end{figure}

\begin{proposition}
There is a cubic graph $H$ which has no 1-factor and $4 \in {\cal S}(H)$.
\end{proposition}

{\bf Proof.} The graph $H$ of Figure \ref{best} has a nowhere-zero 4-flow. It has no 1-factor since $H-v$ has more than one odd component. 
By Theorem  \ref{main_theorem_2}, it has no nowhere-zero 3-flow. Hence, $4 \in {\cal S}(G)$ by 
Theorem \ref{gap}.\hfill $\square$

\begin{corollary} \label{k>5}
Let $G$ be a cubic graph that does not have a 1-factor.
If $k \in \overline{{\cal S}}(G)$, then $k \geq 5$. 
\end{corollary}

\subsection*{Smallest $r$-minimal sets}

Let $r \geq 2$. An $r$-minimal  set $X$ is a {\em smallest}
$r$-minimal set of $G$ if $|X| \leq |X'|$ for every $r$-minimal set $X'$ of $G$. 

\begin{proposition} \label{smallest}
Let $t \geq 1$ be an integer, $G$ a $(2t+1)$-regular graph and $r \geq 2$. If $X \subseteq E(G)$ is a smallest $r$-minimal set, 
then $\Delta(G[X]) \leq t$. 
\end{proposition}
{\bf Proof.} Suppose to the contrary that $\Delta(G[X]) > t$. Then there
is $v \in V(G)$ such that $d_{G[X]}(v) > t$. If we switch at $v$, then we obtain an equivalent graph $(G,\sigma)$
with $|N_{\sigma}| < |X|$ and $F_c((G,\sigma)) = r$. But $N_{\sigma}$ contains an $r$-minimal set $X'$, contradicting the fact that 
$X$ is a smallest $r$-minimal set. \hfill $\square$

We will prove some bounds for the cardinality of 
smallest $r$-minimal sets. The independence number of $G$ is denoted by $\alpha(G)$.

\begin{proposition} \label{upper_bound_smallest_set}
Let $t \geq 1$ be an integer and $G$ be a $(2t+1)$-regular graph. If $X \subseteq E(G)$ is a 
smallest $(2 + \frac{1}{t})$-minimal set, then
$|X| \leq \min\{(\frac{1}{2}|V(G)| - \alpha(G))(2t+1), \frac{t}{2}|V(G)|\}$.
\end{proposition}
{\bf Proof.} Since $G$ has a $(2 + \frac{1}{t})$-minimal set, it follows by Theorem \ref{always_2+1/t} that
$G$ has a $t$-factor. 
Let $(G,\sigma)$ be the graph with $N_{\sigma} = E(G)$ and $V \subseteq V(G)$ be an independent set with
$|V| = \alpha(G)$. The function $\phi : E(G) \rightarrow \{1\}$ is a modular nowhere-zero $(2+\frac{1}{t})$-flow on $(G,\sigma)$.
Switch at every vertex of $V$ to obtain an equivalent graph $(G,\sigma')$ with a modular nowhere-zero 
$(2+\frac{1}{t})$-flow $\phi'$ and $|N_{\sigma'}| \leq |E(G)| - \alpha(G)(2t+1) = (\frac{1}{2}|V(G)| - \alpha(G))(2t+1)$. 
It follows with Theorem \ref{gen_Xu_Zhang} that $(G,\sigma')$ has a nowhere-zero $(2 + \frac{1}{t})$-flow. Therefore, 
$|X| \leq (\frac{1}{2}|V(G)| - \alpha(G))(2t+1)$. Proposition \ref{smallest} implies that $|X| \leq \frac{t}{2}|V(G)|$.
\hfill $\square$

Let $G$ be a bridgeless cubic graph. 
The {\em resistance} $r(G)$ is the cardinality of a minimum
color class, where the minimum is taken over all proper 4-edge-colorings of $G$. It is easy to see that $r(G) \leq \omega(G)$
and if $r(G) \not = 0$, then $r(G) \geq 2$ (see \cite{Steffen_1998}). A bridgeless cubic graph which is not 3-edge-colorable is
called a {\em snark}.

\begin{theorem} \label{cubic_details}
Let $G$ be a cubic graph which has a 1-factor and $G \not = K_2^3$. For each $i \in \{3,4\}$ there is a smallest $i$-minimal set $X_i$ in $G$, and  \\
1) if $G$ is bipartite, then $|X_3| = 0$ and $|X_4| = 2$. \\
2) if $G$ is 3-edge-colorable and not bipartite, then $ 2 \leq |X_3| \leq 3(\frac{1}{2}|V(G)| - \alpha(G))$ and  $|X_4| = 0$.\\
3.1) if $G$ is not 3-edge-colorable, then $r(G) \leq |X_4| \leq \omega(G) \leq  |X_3| \leq \min \{3(\frac{1}{2}|V(G)| - \alpha(G)), \frac{1}{2}|V(G)|\}$. \\
3.2) if $G$ is a snark, then $r(G) \leq |X_4| \leq \omega(G) <  |X_3| \leq \min \{3(\frac{1}{2}|V(G)| - \alpha(G)), \frac{1}{2}|V(G)|\}$.
\end{theorem}
{\bf Proof.} It follows from Theorem \ref{main_theorem} that there is a smallest $i$-minimal set $X_i$
in $G$ for each $i \in \{3,4\}$. By Proposition \ref{upper_bound_smallest_set},  $|X_3| \leq 3(\frac{1}{2}|V(G)| - \alpha(G))$.

1) If $G$ is bipartite, then $|X_3| = 0$ and it follows with  Theorem \ref{facts_4_minimal_bipartite} that
$G$ has a signature $\sigma$ with $|N_{\sigma}|=2$ and $F((G,\sigma)) = 4$. Since every signature of a flow-admissible graph has
at least two edges it follows that $|X_4| = 2$.

2) If $G$ is not bipartite and 3-edge-colorable, then $|X_4| = 0$, and as above we get that $|X_3| \geq 2$.

3) Let $G$ be not 3-edge-colorable. Suppose to the contrary
that there is a smallest 4-minimal set $X$ such that $|X| < r(G)$.
By Lemma \ref{induced_flow}, $(G_{X}^2,\emptyset)$ has a nowhere-zero 4-flow and hence, it is 3-edge-colorable. 
By Proposition \ref{smallest}, $X$ is an independent set. Hence,
a 3-edge-coloring of $G_{X}^2$ induces a proper 4-edge-coloring of $G$ which has a minimal color class 
$c$ with $|c| \leq |N_{\sigma}| < r(G)$, contradiction. Thus, $r(G) \leq |X_4|$.
The other inequalities follow by Theorem \ref{facts_4_minimal_class2}. \nolinebreak \hfill $\square$

The bounds of Theorem \ref{cubic_details} are sharp for the Petersen graph. 
In \cite{Steffen_2004} it is shown that for every positive integer $k$ there is a cubic graph $G$
such that $\omega(G) - r(G) \geq k$ and that there is cyclically 5-edge-connected cubic graph $H$
and $r(H) \geq k$.
If $G$ is not 3-edge-colorable and $r(G) > 2$, then there is a signature $\sigma$ such that $(G,\sigma)$ is flow-admissible
and $2 \leq |N_{\sigma}| < r(G)$. Theorem \ref{cubic_details} implies $F((G,\sigma)) > 4$. Hence, we obtain the following corollary
which is similar to Corollary \ref{k>5}.

\begin{corollary} Let $G$ be a cubic graph. If $\omega(G) > 2$, then for every $k$ with $2 \leq k < r(G)$ there is a 
signature $\sigma$ such that $|N_{\sigma}| = k$ such that $F_c((G,\sigma)) > 4$ and $F((G,\sigma)) \geq 5$.
\end{corollary}

\subsection*{The integer flow spectrum of a class of cubic graphs}

It is quite difficult to determine the flow spectrum of a graph. Indeed even for the integer flow spectrum it is 
difficult. So far the integer flow spectrum has been determined only for eulerian graphs in \cite{Skoviera_2011}, and for
complete and complete bipartite graphs in \cite{Macajova_Rollova_2011}.
For instance, it is known that $\overline{{\cal S}}(G) = \{3,4,5,6\}$, if $G$ is the Petersen graph. 

For $n \geq 1$, let $G_n$ be the cubic graph which is obtained from a circuit of length $2n$, where every second edge is replaced by
two parallel edges. 

\begin{theorem} \label{6_flow_family}
If $n=1$, then $\overline{{\cal S}}(G_n) = \{3\}$. If $n= 2$, then $\overline{{\cal S}}(G_n) = \{3,4\}$, and if $n \geq 3$, 
then $\overline{{\cal S}}(G_n) = \{3,4,6\}$. 
\end{theorem}
{\bf Proof.} If $n=1$, then $G_1 = K_2^3$. Hence,  $\overline{{\cal S}} (G_1) = \{3\}$.

Let $n \geq 2$. Theorem \ref{main_theorem} implies that  $\{3,4\} \subseteq \overline{{\cal S}}(G_n)$.
Let $V(G_n) = \{ v_0, \dots, v_{2n-1}\}$. For $ i \in \{0, \dots, n-1\}$ the vertices $v_{2i}, v_{2i + 1}$ are connected by two parallel edges and the
vertices $v_{2i + 1}, v_{2i+2}$ are connected by a simple edge (indices are added modulo $2n$). 
Every signature of $G_n$ is equivalent to a signature $\sigma$ where for each $ i \in \{0, \dots, n-1\}$ at most one edge between $v_{2i}, v_{2i + 1}$
is negative and all other edges are positive. We call $\sigma$ a normal signature of $G_n$. We say that $\sigma$ is odd or even, depending on whether
$|N_{\sigma}|$ is odd  or even. 
Hence, we have only to consider the two cases whether $\sigma$ is even or odd. 

If $\sigma$ is even, then $G_n$ is the
union of two balanced eulerian graphs and hence, $F((G,\sigma)) \leq 4$ by Lemma \ref{union_eulerian_graph}. Hence, 
$\overline{{\cal S}}(G_2) = \{3,4\}$. 

It remains to consider the case when $n \geq 3$ and $\sigma$ is odd.  Then $|N_{\sigma}| \geq 3$.
Let $e_1,e_2 \in N_{\sigma}$. There are hamiltonian circuits $C_1$, $C_2$ such that $E(C_i) \cap N_{\sigma} = N_{\sigma} - \{e_i\}$.
Both circuits are balanced and hence, there are nowhere-zero 2-flows 
$\phi_i$ on $C_i$. Let $e_1'$ be the positive edges which is parallel to $e_1$. 
Then $\psi = 2 \phi_1 + \phi_2$ is a 4-flow on $G_n$, and 
$\psi(e) \not = 0$ if $e \in N_{\sigma}$, $\psi(e) \in \{0,1,3\}$ if $e \in E(G_n) - (N_{\sigma} \cup \{e_1'\})$, and $\psi(e_1') = 2$.
Let $\tau_{\psi}$ be the underlying orientation of $H(G)$ for $\psi$.
There is a hamiltonian circuit $C$ that consists only of positive edges, that contains $e_1'$ and all edges $e$ with $\psi(e)=0$.
Let $\tau$ be an orientation of $H(C)$ such that $C$ is a directed circuit and $\tau$ and $\tau_{\psi}$ coincide on $e_1'$. 
Let $\psi'$ be a nowhere-zero 2-flow on $C$ with orientation $\tau$. Then,
$\psi + 2 \psi'$ is a nowhere-zero 6-flow on $G_n$.  

Suppose to the contrary that $F((G_n,\sigma)) < 6$. Let $k < 6$ and $\psi_n$ be nowhere-zero $k$-flow on $G_n$.
Without loss of generality we assume in the following that all flow values are positive.

We first show that if $n > 3$, then $F((G_n,\sigma)) < 6$ implies that there is a $m < n$ such that $F((G_m,\sigma_m)) < 6$, 
where $\sigma_m$ is normal and odd.  

Let $n > 3$. It is easy to see that if one edge of the two edges between $v_{2i}$ and $v_{2i+1}$ is negative, then
$\psi_n(v_{2i-1}v_{2i}) \not = \psi_n(v_{2i+1}v_{2i+2})$. 
If there are $i,j \in \{0, \dots, n-1\}$ with $i < j$ and $\psi_n(v_{2i+1}v_{2i+2}) = \psi_n(v_{2j+1}v_{2j+2})$,
then remove these two edges and add edges $v_{2i+2}v_{2j+1}$ and $v_{2i+1}v_{2j+2}$ to obtain two graphs $G_{n_1}$
and $G_{n_2}$ with nowhere-zero $k$-flows $\psi_{n_1}$ and $\psi_{n_2}$, respectively. Depending on the orientation 
of the half-edges of  $v_{2i+1}v_{2i+2}$ and of $v_{2j+1}v_{2j+2}$ the new edges might be negative.  For one of these two graphs, 
say $G_{n_1}$, $\psi_{n_1}$ is equivalent to a $k$-flow $\psi_{n_1}'$ on an odd normal signature $\sigma_{n_1}'$ of $G_{n_1}$, 
since for otherwise $G_n$ and an even normal signature which is equivalent to $\sigma_{n}$ could be reconstructed 
from $(G_{n_1},\sigma_{n_1})$ and $(G_{n_2},\sigma_{n_2})$. 
Hence, $|N_{\sigma_{n_1}}| \geq 3$. 
Thus, we can assume that $G_n$ has an odd normal signature $\sigma_n$ with $|N_{\sigma}| = n$. In particular, $n$ is odd. 

Since $k < 6$, it follows that if $G_n$ is not reducible to a smaller graph $G_m$, then $n=3$.  
Consider $(G_3,\sigma_3)$. Since $\sigma_3$ is normal it follows that the difference between the flow values of any two simple edges 
is at least 2. Since $\psi_3(e) \geq 1$ for every 
edge $e \in G_3$, it follows that there is a simple edge with $\psi_3(e) \geq 5$, contradicting our assumption that $k < 6$.

Therefore $F((G_n,\sigma)) = 6$ and $\overline{{\cal S}}(G_n) = \{3,4,6\}$ if $n > 2$.
\hfill $\square$

The graphs $G_n$ are bipartite. Hence, the empty set is a smallest 3-minimal set for all $n \geq 1$. 
If  $n \geq 2$, then a smallest 4-minimal set contains precisely two edges, and if $n \geq 3$, then a smallest 6-minimal set consists of
three edges.

\subsection*{Bouchet's conjecture}

To prove Bouchet's conjecture it suffices to prove it for cubic graphs. Hence, we are in a similar situation as for Tutte's 5-flow conjecture.
However, as Section \ref{difference} shows Bouchet's conjecture has to be
proved for integer flows explicitly. 

\begin{lemma}\label{union_eulerian_graph}
Let $(G,\sigma)$ be a signed graph. If $(G,\sigma)$ is the union of two balanced eulerian graphs, then $F((G,\sigma)) \leq 4$.
\end{lemma}
{\bf Proof.} Let $H_1$ and $H_2$ be eulerian graphs such that $E(H_1) \cup E(H_2) = E(G)$. For $i \in \{1,2\}$, let $\sigma_i$ be the restriction 
of $\sigma$ to $H_i$. Since $H_i$ is balanced it follows that
there is a nowhere-zero 2-flow $\phi_i$ on $(H_i,\sigma_i)$. Hence, $\phi_1 + 2\phi_2$ is a nowhere-zero 4-flow on $(G,\sigma)$.
\hfill $\square$

\begin{theorem} \label{6-flow_1}
Let $(G,\sigma)$ be a flow-admissible signed cubic graph. If $G$ is a Kotzig-graph, then $F((G,\sigma)) \leq 6$.
\end{theorem}
{\bf Proof.} Since $G$ is a Kotzig-graph, $G$ has three 1-factors $M_1$, $M_2$, $M_3$ such that the union of any two of them
induces a hamiltonian circuit of $G$. It follows that there are two, say $M_1$ and $M_2$ such that 
$|N_\sigma \cap M_1|$ and $|N_\sigma \cap M_2|$ have the same parity. Hence, $|N_\sigma \cap M_1| + |N_\sigma \cap M_2|$ is even,
and $(G[M_1 \cup M_2],\sigma_{1,2})$ with $\sigma_{1,2} = \sigma|_{M_1 \cup M_2}$ is balanced. Clearly, $(G,\sigma)$ is equivalent to
$(G,\sigma')$ with $N_{\sigma'} \cap (M_1 \cup M_2) = \emptyset$.

If $|M_3 \cap N_{\sigma'}|$ is even, then $(G[M_1 \cup M_3],\sigma_{1,3})$ with 
signature $\sigma_{1,3} = \sigma'|_{M_1 \cup M_3}$ is balanced. Hence, $F(G,\sigma) \leq 4$ by Lemma \ref{union_eulerian_graph}.

If $|M_3 \cap N_{\sigma'}|$ is odd, then $|N_{\sigma'}| \geq 3$. Let $e=xy$ be an extroverted edge of $(G,\sigma')$, 
and let $(G,\sigma^*)$ be the graph which is obtained from $(G,\sigma')$ by changing the direction of $h_e^x$. Then 
$e$ is a positive edge which is directed from $x$ to $y$ in $(G,\sigma^*)$.
It follows as above, that $(G,\sigma^*)$ has a nowhere-zero 4-flow $\phi$. Without loss of generality we can assume 
that all flow values are positive, all edges of $M_3$ have flow value 1, and $M_1 \cup M_2$ is a directed circuit $C$.
For two vertices $a,b$, let $P(a,b)$ denote the directed path from $a$ to $b$ in $C$.

If we consider $\phi$ on $(G,\sigma')$, then $\delta\phi(x)=2$, and $\delta\phi(v)=0$ for all $v \in V(G) \setminus \{x\}$. 
Since all flow values are positive, it follows that there is an introverted edge $f = uw$. Let $x,u,w$ be the sequent order of these three vertices 
in $C$. 
We define a nowhere-zero 6-flow $\phi^*$ on $(G,\sigma')$ as follows:
$\phi^*(e) = \phi(e)$, if $e \in E(G) - (E(P(x,w)) \cup \{f\})$, $\phi^*(e) = \phi(e) + 2$, if $e \in E(P(x,u))$,
$\phi^*(e) = \phi(e) + 1$, if $e \in E(P(u,w))$, and  $\phi^*(f) = 2$.  \hfill $\square$

Greenwell and Kronk \cite{Greenwell_1973} proved that every uniquely 3-edge-colorable cubic graph has precisely three hamiltonian
circuits, which are induced by the three color classes. Hence, we obtain the following corollary.

\begin{corollary}
Let $(G,\sigma)$ be a flow-admissible signed cubic graph. If $G$ is uniquely 3-edge-colorable, then $F((G,\sigma))\leq 6$.
\end{corollary}

Further results and references on Kotzig-graphs and uniquely 3-edge-colorable graphs can be found in \cite{Zhang_2012}.

Remark: For $t \geq 2$ it can be proved by a slight modification of the proof of Theorem \ref{6-flow_1} that $F((G,\sigma)) \leq 4$ for every 
flow-admissible $(2t+1)$-regular Kotzig-graph $(G,\sigma)$. However, this is also a simple consequence of a result of Raspaud and Zhu 
\cite{Raspaud_Zhu_2011} that every flow-admissible 4-edge-connected graph has a nowhere-zero 4-flow.

\section{Concluding remarks}

\begin{figure}
\centering
	\includegraphics[height=5cm]{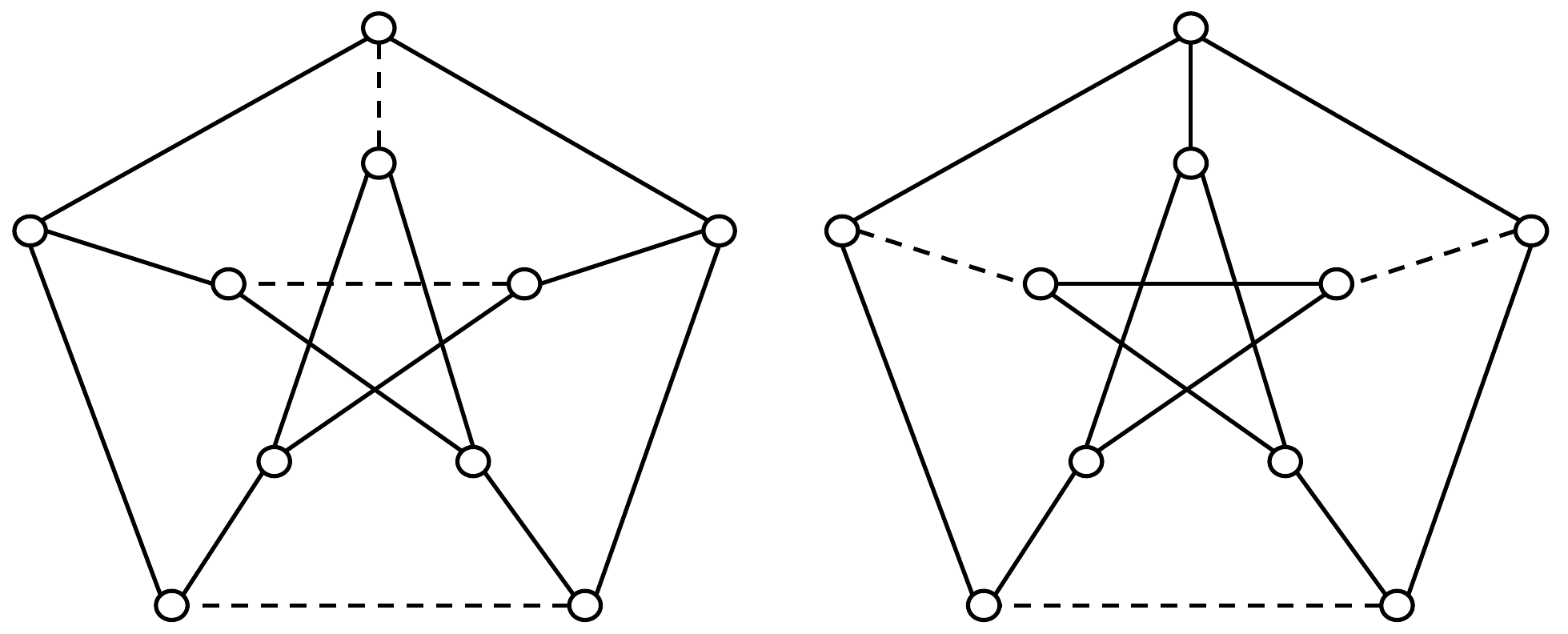}
\caption{the dotted lines represent set $X$ on the left and $X'$ on the right\label{P}} \end{figure}

If a graph $H$ has an $r$-minimal set $X$ of  cardinality 2, then ${\cal S}_X(H) = \{r, F_c((H,\emptyset))\}$. 
Theorem \ref{cubic_details} implies that every 3-minimal set of a bridgeless non-3-edge-colorable cubic graph contains at
least three edges. The Petersen graph $P$ has a 3-minimal set $X$ with $|X| = 3$. Hence, ${\cal S}_X(P) = \{3,4,5\}$.
On the other hand, $P$ has a 6-minimal set $X'$ with $|X'| = 3$. It follows that $3 \not \in {\cal S}_{X'}(P)$; 
indeed $\overline{{\cal S}}_{X'}(P) = \{4,5,6\}$. 
The sets $X$ and $X'$ are indicated in Figure \ref{P}. 
Two switches at $v$ and $w$ yield a signature $\sigma$ of $P$ such that $(P,\sigma)$ and $(P,X)$ are equivalent. 
But $X' \subset N_{\sigma}$ and therefore, ${\cal S}_{X}(P) \not = {\cal S}_{N_{\sigma}}(P)$.
We conclude with the following questions.

\begin{problem}
Let $r \geq 2$, $(G,\sigma)$ be a flow-admissible signed graph and $X$ a (non-empty) $r$-minimal set. 
Determine the (integer) $X$-flow spectrum of $G$.
\end{problem}

\begin{problem}
Let $r \geq 2$ and $G$ be a graph. Let $(G,\sigma)$ and $(G,\sigma')$ be flow-admissible. Is it true that
if $N_{\sigma}$ and $N_{\sigma'}$ are (both) smallest $r$-minimal sets, then ${\cal S}_{N_{\sigma}}(G) = {\cal S}_{N_{\sigma'}}(G)$ 
$($or $\overline{{\cal S}}_{N_{\sigma}}(G) = \overline{{\cal S}}_{N_{\sigma'}}(G)$, if $r$ is an integer$)$?
\end{problem}

Theorem \ref{cubic_details} says that for snarks every 3-minimal set contains a proper 4-minimal set. If Tutte's 5-flow conjecture
is true, then the following problem for integer flows has an affirmative answer. 

\begin{problem} Let $G$ be a snark and $X_3$ a 3-minimal set. Is it true that there exist a 4-minimal set $X_4$ and 
a 5-minimal set $X_5$ such that $X_5 \subset X_4 \subset X_3$?
\end{problem}


\begin{thebibliography}{99}


\bibitem{Akiyama}{\sc J.~Akiyama, M.~Kano}, Factors and factorizations of graphs, Lecture Notes in Mathematics 2031, Springer-Verlag Berlin, 
		Heidelberg (2011)

\bibitem{Bouchet} {\sc A.~Bouchet}, Nowhere-zero integral flows on bidirected graph, J.~Comb. Theory Ser.~B {\bf 34} (1983) 279 - 292


\bibitem{Greenwell_1973} {\sc D.~Greenwell, H.V.~Kronk}, Uniquely line-colorable graphs, Can.~Math. Bull.{\bf 16} (1973) 525 - 529

\bibitem{Macajova_Rollova_2011} {\sc E.~M\'a\v{c}ajov\'a, E.~Rollov\'a}, On the flow numbers of signed complete and complete bipartite graphs,
	Elect.~Notes in Disc.~Math.~{\bf 38} (2011) 591-596

\bibitem{Skoviera_2011} {\sc E.~M\'a\v{c}ajov\'a, M.~\v{S}koviera}, Nowhere-zero flows on signed eulerian graphs, Technical Reports in Informatics
	TR-2011-029, Faculty of Mathematics, Physics, and Informatics, Comenius University, Bratislava (2011)

\bibitem{Petersen_1891} {\sc J.~Petersen}, Die Theorie der regul\"aren graphs, Acta Mathematica {\bf 15} (1891) 193 - 220  

\bibitem{Raspaud_Zhu_2011} {\sc A.~Raspaud, X.~Zhu}, Circular flow on signed graphs, J.~Comb.~ Theory Ser.~B {\bf 101} (2011) 464 - 479

\bibitem{Seymour_6-flow} {\sc P.D.~Seymour}, Nowhere-zero 6-flows, J.~Comb.~Theory Ser.~B {\bf 30} (1981) 130 - 135

\bibitem{Steffen_1998} {\sc E.~Steffen}, Classifications and characterizations of snarks, Disc. Math.~{\bf 188} (1998) 183 - 203

\bibitem{Steffen_2001}{\sc E.~Steffen}, Circular flow numbers of regular multigraphs, J. Graph Theory {\bf 36} (2001) 24 - 34 

\bibitem{Steffen_2004} {\sc E.~Steffen}, Measurements of edge-uncolorability, Disc.~Math.~{\bf 280} (2004) 191 - 214

\bibitem{Tutte_1949} {\sc W.T.~Tutte}, On the immbedding of linear graphs in surfaces, Proc.~London Math.~Soc.~{\bf 22} (1949) 474 - 483

\bibitem{Tutte_1954} {\sc W.T.~Tutte}, A contribution to the theory of chromatic polynomials, Canad. J. Math. {\bf 6} (1954) 80 - 91

\bibitem{XuZhang} {\sc R.~Xu, C.-Q.~Zhang}, On flows in bidirected graphs, Disc.~Math.~{\bf 299} (2005) 335 - 343

\bibitem{Zhang_2012} {\sc C.-Q.~Zhang}, Circuit Double Cover of Graphs,  London Mathematical Society Lecture Note Series {\bf 399},
	Cambridge University Press, Cambridge (2012) 

\bibitem{Zyka_1987} {\sc O.~Z\'{y}ka}, Nowhere-zero 30-flow on bidirected graphs, Thesis, Charles University, Praha, KAM-DIMATIA Series 87-26 (1987)

\end{thebibliography}
\end{document}